\documentclass[11pt]{extarticle}
\usepackage[T1]{fontenc}
\usepackage{amsmath}
\usepackage{amsthm}
\usepackage{amsfonts}
\usepackage{amssymb}
\usepackage{xcolor}
\usepackage{graphicx}
\usepackage{caption}
\usepackage{multirow}
\usepackage{placeins}

\usepackage[mathlines]{lineno}

\numberwithin{equation}{section}
\theoremstyle{plain}
\newtheorem{theorem}{Theorem}[section]
\newtheorem{lemma}[theorem]{Lemma}
\newtheorem{corollary}[theorem]{Corollary}
\newtheorem{proposition}[theorem]{Proposition}
\theoremstyle{remark}
\newtheorem{remark}[theorem]{Remark}

\begin{document}

\title{A Hardy space approximation supporting zero-free half-planes for the $\zeta$-function}

\author{
    Juan Manzur\thanks{Department of Mathematics,
    Universidad del Atlántico, Puerto Colombia, Atlántico, Colombia. email: jcmanzur@mail.uniatlantico.edu.co}
    \and
    Gustavo Quintero\thanks{Department of Mathematics, Universidad del Atlántico, Puerto Colombia, Atlántico, Colombia. email: gdquintero@mail.uniatlantico.edu.co}
}

\date{}
\maketitle

\begin{abstract}
By a theorem of Bagchi, the Riemann Hypothesis (RH) holds if and only if
$E(s) = 1/s$ belongs to the closed linear span of
$G_k(s) = (k^{-s} - k^{-1})\zeta(s)/s$, $k \geq 2$, in the Hardy space
$H^2(\mathbb{C}_{1/2})$, where $\mathbb{C}_\alpha$ is the half-plane
$\mathrm{Re}(s) > \alpha$. The formal M\"obius expansion
$E = \sum_{k \geq 2} \mu(k) G_k$ diverges on the critical line, by a theorem of
B\'aez-Duarte. We study it in the scale of spaces $H^2(\mathbb{C}_\alpha)$, where its
convergence forces $\zeta$ to be zero-free in $\mathbb{C}_\alpha$. We prove that the
partial sums $S_n$ converge to $E$ in $H^2(\mathbb{C}_\alpha)$ for every $\alpha > 1$,
and that a single zero of $\zeta$ in $\mathbb{C}_\alpha$ keeps the best-approximation
error $d_n(\alpha)^2$ uniformly away from zero. We express the M\"obius error
$s_n(\alpha) = \|S_n - E\|^2_\alpha$ as a series with elementary terms, free of
$\zeta$, and compute it for $n \leq 10^6$: $s_n(1/2)$ oscillates without decaying,
whereas for $\alpha > 1/2$ it decays at a rate close to $n^{1-2\alpha}$, much faster
than the bound known under RH. On the critical line, $d_n(1/2)^2$ behaves like
$0.046/\log n$ for $100 \leq n \leq 5120$, in agreement with earlier computations by
B\'aez-Duarte, Balazard, Landreau and Saias.

\bigskip
\noindent
\textbf{Keywords:} Riemann hypothesis, Nyman--Beurling criterion, zero-free half-planes, Hardy spaces, numerical approximation.

\medskip
\noindent
\textbf{Mathematics Subject Classification:} 11M26, 30H10, 46E22, 11Y35.
\end{abstract}

\section{Introduction}\label{sec:intro}

The Riemann Hypothesis (RH) asserts that all non-trivial zeros of the Riemann
zeta function lie on the critical line $\mathrm{Re}(s) = \frac{1}{2}$. Its
implications for the distribution of prime numbers make it one of the central
open problems in mathematics.

A functional-analytic approach to RH was initiated by Nyman and
Beurling~\cite{beurling,Nyman}. Writing $\{y\}$ for the fractional part of $y$
and $\chi = \chi_{(0,1]}$, their criterion states that RH holds if and only if
$\chi$ can be approximated in $L^2(0,\infty)$ by linear combinations of the
dilations $x \mapsto \{\theta/x\}$, $0 < \theta < 1$; see
\cite[Theorem~1.1 and Remark~1.1]{baez-unitary} for this form.
Báez-Duarte~\cite{baez-beurling,baez-unitary} examined the most natural candidate
for such an approximation. M\"obius inversion of the identity
$\lfloor 1/x \rfloor = \sum_{j \geq 1} \chi(jx)$ \cite[(2.4)]{baez-unitary} gives
$\chi(x) = \sum_{j \geq 1} \mu(j) \lfloor 1/(jx) \rfloor$
\cite[(2.6)]{baez-unitary}, where $\mu$ is the
M\"obius function, and writing $\lfloor y \rfloor = y - \{y\}$ and using
$\sum_{j} \mu(j)/j = 0$ one arrives at
\begin{equation}\label{eq:bd-series}
\chi(x) = -\sum_{j=1}^{\infty} \mu(j) \Bigl\{ \frac{1}{jx} \Bigr\}.
\end{equation}
This series converges pointwise and in $L^1(0,1)$~\cite[\S2.2]{baez-unitary}, but
Báez-Duarte proved that it
does not converge in $L^2(0,\infty)$~\cite[Theorem~2.2]{baez-unitary}: the
M\"obius coefficients, natural as they are, fail precisely in the norm in which
the criterion is formulated.

The same coefficients run through the subsequent literature. Báez-Duarte later
strengthened the Nyman--Beurling criterion by replacing the uncountable family of
dilations with the countable family $x \mapsto \{1/(kx)\} - \{1/x\}/k$,
$k \geq 2$~\cite{baez}. In an
equivalent weighted $\ell^2$-formulation, in which the constant sequence
$\mathbf{1}$ is approximated by the remainder sequences $r_k(j) = j \bmod k$,
Balazard and Saias~\cite{Balazard} asked whether the coefficients $c_{n,k}$ of the
orthogonal projection of $\mathbf{1}$ onto
$\mathrm{span}\{r_k : 2 \leq k \leq n\}$ satisfy
$\lim_{n\to\infty} c_{n,k} = -\mu(k)/k$ for every $k \geq 2$, and
Weingartner~\cite{Weingartner} proved that this is equivalent to RH. In the Hardy
space $H^2(\mathbb{D})$ of the unit disk, Noor~\cite{waleed1} showed that RH holds
if and only if the constant function $\mathbf{1}$ belongs to the closed linear
span of
\[
h_k(z) = \frac{1}{1-z}\log\!\left(\frac{1+z+\cdots+z^{k-1}}{k}\right),
\quad k \geq 2,
\]
and proved that $\sum_{k \geq 2} \frac{\mu(k)}{k}(I-S)h_k$ converges to $1-z$,
where $S$ is the shift operator. Ghosh, Kremnizer, Noor and
Santos~\cite{noor-charles} developed a general approach in which the density of
$\mathrm{span}\{h_k\}$ in suitable spaces of analytic functions yields zero-free
half-planes, and showed in particular that a M\"obius series in the $h_k$ converges
in $H^q(\mathbb{D})$, $0 < q < 1$, which gives the half-plane $\mathrm{Re}(s) > 1$;
see also~\cite[pp.~5--6]{noor-charles} for an overview of the many equivalent
reformulations of the Nyman--Beurling and B\'aez-Duarte criteria. Throughout, the M\"obius coefficients are the
canonical candidates, and what is at stake is the norm in which convergence is
measured.

In this paper we pursue this question in the Hardy spaces of half-planes. For
$\alpha \in \mathbb{R}$ let
$\mathbb{C}_\alpha = \{s \in \mathbb{C} : \mathrm{Re}(s) > \alpha\}$. Bagchi
\cite{bagchi} showed that RH holds if and only if $E(s) = 1/s$ belongs to the
closed linear span of
\[
G_k(s) = \bigl(k^{-s} - k^{-1}\bigr)\frac{\zeta(s)}{s}, \quad k \geq 2,
\]
in $H^2(\mathbb{C}_{1/2})$. The Mellin transform carries $\chi$ to $E$ and
$\{1/(kx)\} - \{1/x\}/k$ to $-G_k$ (see Subsection~\ref{subsec:method}), and under
this correspondence~\eqref{eq:bd-series} becomes the formal identity
\begin{equation}\label{eq:formal}
E = \sum_{k=2}^{\infty} \mu(k)\, G_k ,
\end{equation}
which also follows formally from $\sum_{k} \mu(k) k^{-s} = 1/\zeta(s)$ and
$\sum_{k} \mu(k)/k = 0$. By Plancherel's theorem the norm of
$H^2(\mathbb{C}_{1/2})$ corresponds to the norm of $L^2(0,1)$, and Báez-Duarte's
theorem shows that the partial sums $S_n = \sum_{k=2}^{n} \mu(k) G_k$ do not
converge to $E$ in $H^2(\mathbb{C}_{1/2})$ (Remark~\ref{rem:exact}). The space
$H^2(\mathbb{C}_{1/2})$, however, is only the first of the scale
$H^2(\mathbb{C}_\alpha)$, and since $\|F\|_\alpha \leq \|F\|_{1/2}$ for
$\alpha \geq 1/2$, where $\|\cdot\|_\alpha$ denotes the norm of
$H^2(\mathbb{C}_\alpha)$, convergence becomes easier as $\alpha$ increases. This
raises the question that motivates the present work: for which $\alpha$ does
\eqref{eq:formal} converge in $H^2(\mathbb{C}_\alpha)$, and what does its
convergence tell us about $\zeta$?

The starting point is the following simple observation, which ties convergence in
a half-plane to the zeros of $\zeta$ in that half-plane
(Proposition~\ref{main-prop}): for every $\alpha \geq 1/2$,
\[
E \in \overline{\mathrm{span}}\{G_k : k \geq 2\} \ \text{in }
H^2(\mathbb{C}_\alpha) \implies \zeta(s) \neq 0, \quad \forall\, s \in \mathbb{C}_\alpha .
\]
In particular, the convergence of the M\"obius series~\eqref{eq:formal} in
$H^2(\mathbb{C}_\alpha)$ is a sufficient condition for the half-plane
$\mathbb{C}_\alpha$ to be free of zeros, and Báez-Duarte's obstruction is the case
$\alpha = 1/2$ of a whole scale of such conditions. We then prove
(Theorem~\ref{main-teo}) that $S_n \to E$ in $H^2(\mathbb{C}_\alpha)$ for every
$\alpha > 1$; in the scale $H^2(\mathbb{C}_\alpha)$ this is the counterpart of the
result of~\cite{noor-charles} for $H^q(\mathbb{D})$. The proof extends Noor's
convergence result to a family of functions $g_{k,\tau}$ in $H^2(\mathbb{D})$ and
transfers it, through the Hilbert space $\mathcal{H}^2$ of Dirichlet series, to
$H^2(\mathbb{C}_\alpha)$ by a bounded operator
$\mathcal{M}^\tau : \mathcal{H}^2 \to H^2(\mathbb{C}_{1/2+\tau})$. A
quantitative counterpart also holds. For $n \geq 2$ and $\alpha > 0$ let
\[
d_n(\alpha)^2 := \inf_{c_2, \dots, c_n \in \mathbb{C}}
\Bigl\| E - \sum_{k=2}^{n} c_k G_k \Bigr\|^2_\alpha
\]
denote the optimal approximation error. Theorem~\ref{thm:lower} shows that a single
zero of $\zeta$ in $\mathbb{C}_\alpha$ keeps $d_n(\alpha)^2$ bounded away from $0$,
uniformly in $n$, and Corollary~\ref{cor:sensor} shows that the threshold
$\alpha_c = \inf\{\alpha \geq 1/2 : d_n(\alpha)^2 \to 0\}$ satisfies
$1/2 \leq \Theta \leq \alpha_c \leq 1$, where $\Theta$ is the supremum of the real
parts of the non-trivial zeros of $\zeta$, with $\alpha_c = 1/2$ if and only if RH holds.

Between the critical line, where~\eqref{eq:formal} diverges, and the half-plane
$\mathrm{Re}(s) > 1$, where it converges, lies the critical strip, and it is there
that the question becomes decisive. By Proposition~\ref{main-prop}, a proof of
convergence for some $\alpha < 1$ would yield a zero-free half-plane far beyond the
known zero-free regions of $\zeta$, so we explore the strip numerically. An exact
representation of the M\"obius error
\[
s_n(\alpha) = \|S_n - E\|^2_\alpha
\]
as a series with elementary terms (Proposition~\ref{prop:exact}), which avoids both
the zeta function and numerical quadrature, allows us to compute $s_n(\alpha)$ up to
$n = 10^6$ for $1/2 \leq \alpha \leq 1$. On the critical line $s_n(1/2)$ oscillates
without decaying, consistently with Báez-Duarte's theorem, which implies that
$s_n(1/2) \not\to 0$ (Remark~\ref{rem:exact}). For every $\alpha > 1/2$, on the
other hand, our computations show that the oscillations decay, at a rate close to
$n^{1-2\alpha}$ over the computed range. Under RH, the convergence
$s_n(\alpha) \to 0$ for $\alpha > 1/2$ follows from the estimates of Balazard and
Saias~\cite[Lemmes~4 et~5]{Balazard}, but only with a bound of order
$n^{-(\alpha-1/2)/3}$; the computations thus indicate a much faster decay, matching
the order of the lower bound in Remark~\ref{rem:exact}(ii). Finally,
in the spirit of Weingartner's theorem, which compares the optimal coefficients
with the M\"obius ones in the $\ell^2$-setting recalled above, we compare on the
critical line the
M\"obius error $s_n(1/2)$ with the optimal error $d_n(1/2)^2$, the distance from $E$
to $\mathrm{span}\{G_2, \dots, G_n\}$. While $s_n(1/2)$ keeps oscillating,
$d_n(1/2)^2$ behaves like $0.046/\log n$ for $100 \leq n \leq 5120$, in agreement with
the computations of B\'aez-Duarte, Balazard, Landreau and Saias~\cite{BBLS}. Since
$d_n(1/2)^2 \to 0$ under RH, whereas $s_n(1/2) \not\to 0$ unconditionally, under RH
the failure at $\alpha = 1/2$ lies in the fixed choice of coefficients $\mu(k)$
rather than in the subspaces themselves. This is in line with the work of Bettin,
Conrey and Farmer~\cite{BCF}, who showed, in the closely related setting
of~\cite{BBLS}, that under RH and a hypothesis on $\sum_\rho |\zeta'(\rho)|^{-2}$ the
smoothed M\"obius coefficients $\mu(k)(1 - \log k/\log n)$ attain the conjectured
optimal rate $(2 + \gamma - \log 4\pi)/\log n$.

The main contributions of the paper are thus: the convergence of the M\"obius
expansion in $H^2(\mathbb{C}_\alpha)$ for every $\alpha > 1$, obtained by an
operator-theoretic transfer of Noor's convergence result; an exact representation
of the M\"obius error, free of $\zeta$, valid for every $\alpha > 0$; and numerical
evidence, up to $n = 10^6$, that for $\alpha > 1/2$ this error decays like
$n^{1-2\alpha}$, much faster than the known conditional bound.

Section~\ref{sec:prelim} introduces the function spaces and the results needed
later. Section~\ref{sec:approx} develops the theoretical framework and proves the
results described above. Section~\ref{sec:numerics} presents the numerical
experiments and their implications.

\section{Preliminaries}\label{sec:prelim}

We collect the function spaces and key results needed throughout the paper. The three Hilbert spaces $H^2(\mathbb{D})$, $\mathcal{H}^2$, and $H^2(\mathbb{C}_{1/2})$ are introduced here as the analytic framework within which the completeness conditions equivalent to the Riemann Hypothesis are formulated and studied.

\subsection{Hardy--Hilbert spaces and Dirichlet series}

Let $H^2(\mathbb{D})$ be the Hardy--Hilbert space of the unit disk, defined as
\[
H^2(\mathbb{D}) = \left\{ f(z) = \sum_{n=0}^{\infty} a_n z^n : \|f\|_{H^2}^2 = \sum_{n=0}^{\infty} |a_n|^2 < \infty \right\}.
\]
A recent reformulation of the RH, introduced in~\cite{waleed1}, is stated as follows. For each $k \geq 2$, define
\[
h_k(z) = \frac{1}{1-z} \log\!\left(\frac{1 + z + \cdots + z^{k-1}}{k}\right).
\]

\begin{theorem}[Theorem 6, \cite{waleed1}]
The Riemann Hypothesis holds if and only if the constant function $\mathbf{1}$ belongs to the closed linear span of $\{h_k : k \geq 2\}$ in $H^2(\mathbb{D})$.
\end{theorem}

Let $S = M_z$ denote the shift operator on $H^2(\mathbb{D})$, corresponding to multiplication by $z$. An additional result from~\cite{waleed1} plays a central role in this paper.

\begin{lemma}[Lemma 11, \cite{waleed1}]
\label{lem1:waleed}
The series $\sum_{k=2}^{\infty} \frac{\mu(k)}{k}(I - S)h_k$ converges to $1 - z$ in $H^2(\mathbb{D})$, where $\mu$ is the M\"{o}bius function.
\end{lemma}

Let $\mathcal{H}^2$ denote the Hilbert space of Dirichlet series~\cite{2}, defined as
\[
\mathcal{H}^2 = \left\{ f(s) = \sum_{n=1}^{\infty} a_n n^{-s} : \|f\|_{\mathcal{H}^2}^2 = \sum_{n=1}^{\infty} |a_n|^2 < \infty \right\}.
\]
Setting $H^2_0(\mathbb{D}) = H^2(\mathbb{D}) \ominus \mathbb{C}$, there exists an isometric isomorphism~\cite{2}
\[
\psi : H^2_0(\mathbb{D}) \longrightarrow \mathcal{H}^2, \qquad z^k \longmapsto k^{-s}.
\]
In what follows we write $H^2_0 = H^2_0(\mathbb{D})$.

\subsection{The Hardy spaces $H^2(\mathbb{C}_\alpha)$}

For $\alpha \in \mathbb{R}$, let $\mathbb{C}_\alpha = \{s = \sigma + it : \sigma > \alpha\}$ and let $H^2(\mathbb{C}_\alpha)$ denote the Hardy space of analytic functions $F$ on $\mathbb{C}_\alpha$ satisfying
\[
\|F\|^2_{H^2(\mathbb{C}_\alpha)} := \sup_{x > \alpha} \frac{1}{2\pi} \int_{-\infty}^{\infty} |F(x + it)|^2 \, dt < \infty.
\]
Any such $F$ admits a non-tangential boundary value $F^*$ almost everywhere on $\{\mathrm{Re}(s) = \alpha\}$, with
\[
\|F\|^2_{H^2(\mathbb{C}_\alpha)} = \frac{1}{2\pi}\int_{-\infty}^{\infty} |F^*(\alpha + it)|^2 \, dt.
\]
Moreover, $H^2(\mathbb{C}_\alpha)$ is a Hilbert space with reproducing kernel
\[
k_w(s) = \frac{1}{s + \overline{w} - 2\alpha}, \quad w \in \mathbb{C}_\alpha;
\]
in particular, convergence in norm implies pointwise convergence in $\mathbb{C}_\alpha$. For brevity we write $\|\cdot\|_\alpha := \|\cdot\|_{H^2(\mathbb{C}_\alpha)}$.

Define the functions $G_k$ and $E$ in $H^2(\mathbb{C}_{1/2})$, for $k \geq 2$, as
\[
G_k(s) = (k^{-s} - k^{-1})\frac{\zeta(s)}{s} \quad \text{and} \quad E(s) = \frac{1}{s}.
\]
The Riemann Hypothesis reformulation for $H^2(\mathbb{C}_{1/2})$ can be stated as follows.

\begin{theorem}[Theorem 2, \cite{bagchi}]\label{thm:bagchi}
RH is true if and only if $E$ belongs to the closed linear span of $\{G_k : k \geq 2\}$ in $H^2(\mathbb{C}_{1/2})$.
\end{theorem}

\section{An $H^2(\mathbb{C}_\alpha)$ approximation approach}\label{sec:approx}

The main goal of this section is to show that $E$ belongs to the closed linear span 
of $\{G_k : k \geq 2\}$ in $H^2(\mathbb{C}_\alpha)$, for every $\alpha > 1$. 
For $\tau > 1/2$ and $k \geq 2$, define
\[
g_{k,\tau}(z) = k^{\tau-1}\sum_{n=1}^{\infty}\frac{z^n}{n^\tau} - \sum_{n=1}^{\infty}\frac{z^{nk}}{n^\tau}.
\]
It is clear that $g_{k,\tau} \in H^2_0$. Furthermore, we observe that
\begin{align*}
(\psi g_{k,\tau})(s) &= k^{\tau-1}\sum_{n=1}^{\infty}\frac{n^{-s}}{n^\tau} - k^{-s}\sum_{n=1}^{\infty}\frac{n^{-s}}{n^\tau}\\
&= (k^{\tau-1} - k^{-s})\zeta(s+\tau)\\
&= k^{\tau-1}f_{k,\tau}(s) \quad \text{for } \mathrm{Re}(s) > 1/2,
\end{align*}
where $f_{k,\tau}(s) = (1 - k^{1-s-\tau})\zeta(s+\tau)$.

\subsection{A bounded operator from $\mathcal{H}^2$ to $H^2(\mathbb{C}_\alpha)$}

The operator $\mathcal{M}_{1/s} : \mathcal{H}^2 \to H^2(\mathbb{C}_{1/2})$, defined by 
$\mathcal{M}_{1/s} : f(s) \mapsto f(s)/s$, is bounded according to~\cite{6}. 
For $\tau \geq 1/2$, we extend this construction by defining 
$\mathcal{M}^\tau : \mathcal{H}^2 \to H^2(\mathbb{C}_{1/2+\tau})$ via $\mathcal{M}^\tau : f \longmapsto f(s - \tau)/s$. 
To see that $\mathcal{M}^\tau$ is well defined and bounded, note that for $\sigma > \tau + 1/2$,
\[
\frac{1}{2\pi}\int_{-\infty}^{\infty}\frac{|f(\sigma-\tau+it)|^2}{|\sigma+it|^2}\,dt
\leq \frac{1}{2\pi}\int_{-\infty}^{\infty}\frac{|f(\sigma-\tau+it)|^2}{|\sigma-\tau+it|^2}\,dt
\leq \|\mathcal{M}_{1/s}f\|^2_{1/2},
\]
and taking the supremum over $\sigma > \tau + 1/2$ yields
\[
\|\mathcal{M}^\tau f\|_{1/2+\tau} \leq \|\mathcal{M}_{1/s}\|\,\|f\|_{\mathcal{H}^2},
\qquad \forall f \in \mathcal{H}^2.
\]
The relevance of $\mathcal{M}^\tau$ to our approach lies in how it acts on the 
functions $1$ and $f_{k,\tau}$. A direct computation gives
\[
(\mathcal{M}^\tau 1)(s) = \frac{1}{s} = E(s),
\]
and
\[
(\mathcal{M}^\tau f_{k,\tau})(s) = \frac{f_{k,\tau}(s-\tau)}{s} 
= (1-k^{1-s})\frac{\zeta(s)}{s} = -kG_k(s).
\]

\subsection{Convergence in $H^2(\mathbb{C}_\alpha)$ for $\alpha > 1$}

We now establish one of the main results of the paper. The key step is the following 
generalization of Lemma~\ref{lem1:waleed} to the family $\{g_{k,\tau}\}$.

\begin{lemma}\label{lem:main}
The series $\sum_{k=2}^{\infty} \frac{\mu(k)}{k^\tau} g_{k,\tau}$ converges 
to $-z$ in $H^2_0$, for each $\tau > 1/2$.
\end{lemma}

\begin{proof}
Since $g_{1,\tau}\equiv 0$, the $k=1$ term may be included, giving
\begin{align*}
\sum_{k=2}^{n} \frac{\mu(k)}{k^\tau} g_{k,\tau}(z)
&= \sum_{k=1}^{n} \frac{\mu(k)}{k^\tau} \left\{ k^{\tau-1}\sum_{j=1}^{\infty} \frac{z^j}{j^\tau} - \sum_{j=1}^{\infty}\frac{z^{jk}}{j^\tau} \right\}\\
&= \sum_{k=1}^{n} \frac{\mu(k)}{k} \sum_{j=1}^{\infty} \frac{z^j}{j^\tau} - \sum_{k=1}^{n} \frac{\mu(k)}{k^\tau} \sum_{j=1}^{\infty}\frac{z^{jk}}{j^\tau}.
\end{align*}
Since $\sum_{k=1}^{\infty} \mu(k)/k = 0$, a statement equivalent to the prime number
theorem~\cite[Ch.~4]{apostol}, the first term vanishes in the limit:
\[
\left\|\sum_{k=1}^{n} \frac{\mu(k)}{k} \sum_{j=1}^{\infty} \frac{z^j}{j^\tau}\right\|_{H^2_0} 
= \left|\sum_{k=1}^{n} \frac{\mu(k)}{k}\right| \left\|\sum_{j=1}^{\infty} \frac{z^j}{j^\tau}\right\|_{H^2_0} \longrightarrow 0.
\]
It therefore suffices to show that
\begin{equation}\label{eq:suffices}
\left\| -\sum_{k=1}^{n} \frac{\mu(k)}{k^\tau} \sum_{j=1}^{\infty}\frac{z^{jk}}{j^\tau} + z \right\|_{H^2_0} \longrightarrow 0.
\end{equation}
A rearrangement of the double sum gives
\begin{align}
\sum_{k=1}^{n} \frac{\mu(k)}{k^\tau} \sum_{j=1}^{\infty}\frac{z^{jk}}{j^\tau}
&= \sum_{j=1}^{n} \frac{z^j}{j^\tau}\left(\sum_{\substack{d|j,\\ 1\leq d\leq n}}\mu(d)\right) + \sum_{j=n+1}^{\infty} \frac{z^j}{j^\tau}\left(\sum_{\substack{d|j,\\ 1\leq d\leq n}}\mu(d)\right)\notag\\
&= \sum_{j=1}^{n} \frac{z^j}{j^\tau}\left\lfloor\frac{1}{j}\right\rfloor + \sum_{j=n+1}^{\infty} \frac{z^j}{j^\tau}\left(\sum_{\substack{d|j,\\ 1\leq d\leq n}}\mu(d)\right), \label{eq:rearrange}
\end{align}
where the last equality uses the identity $\sum_{d|j}\mu(d) = \lfloor 1/j \rfloor$, which
equals $1$ for $j = 1$ and $0$ for $j \geq 2$~\cite[Theorem~2.1]{apostol}. From~\eqref{eq:rearrange},
\[
\sum_{k=1}^{n} \frac{\mu(k)}{k^\tau} \sum_{j=1}^{\infty}\frac{z^{jk}}{j^\tau} = z + \phi_n(z),
\]
where
\[
\phi_n(z) = \sum_{j=n+1}^{\infty} \left(\frac{1}{j^\tau}\sum_{\substack{d|j,\\ 1\leq d\leq n}}\mu(d)\right) z^j.
\]
It remains to show that $\|\phi_n\|_{H^2_0} \to 0$ as $n \to \infty$. 
Letting $\sigma_0(j)$ denote the number of divisors of $j$, we estimate
\[
\|\phi_n\|^2_{H^2_0} = \sum_{j=n+1}^{\infty} \frac{1}{j^{2\tau}}\left|\sum_{\substack{d|j,\\ 1\leq d\leq n}}\mu(d)\right|^2 
\leq \sum_{j=n+1}^{\infty} \frac{1}{j^{2\tau}}\left(\sum_{d|j}1\right)^2
= \sum_{j=n+1}^{\infty} \frac{\sigma_0(j)^2}{j^{2\tau}}.
\]
Since $\sigma_0(j) = o(j^\epsilon)$ for every $\epsilon > 0$, choosing $\epsilon > 0$
such that $\tau > 1/2 + \epsilon$ yields, for some constant $C_\epsilon > 0$,
\[
\|\phi_n\|^2_{H^2_0} \leq C_\epsilon \sum_{j=n+1}^{\infty} \frac{j^{2\epsilon}}{j^{2\tau}} 
= C_\epsilon \sum_{j=n+1}^{\infty} \frac{1}{j^{2\tau - 2\epsilon}} \longrightarrow 0,
\]
since $2\tau - 2\epsilon > 1$.
\end{proof}

For $n \geq 2$ write $S_n := \sum_{k=2}^{n} \mu(k) G_k$ for the $n$-th M\"obius partial sum. With this notation, Lemma~\ref{lem:main}
transfers to the half-plane setting as follows.

\begin{theorem}\label{main-teo}
The sequence $(S_n)$ converges to $E$ in $H^2(\mathbb{C}_{1/2+\tau})$, for
every $\tau > 1/2$.
\end{theorem}

\begin{proof}
By Lemma~\ref{lem:main} and the isomorphism $\psi$,
\[
\sum_{k=2}^{n} \frac{\mu(k)}{k} f_{k,\tau} \longrightarrow -1 \quad \text{in } \mathcal{H}^2.
\]
Applying $\mathcal{M}^\tau$ and using $(\mathcal{M}^\tau f_{k,\tau})(s) = -kG_k(s)$, we conclude
\begin{equation}\label{eq:main-conv}
S_n = \sum_{k=2}^{n} \mu(k) G_k \longrightarrow E \quad \text{in } H^2(\mathbb{C}_{1/2+\tau}).
\end{equation}
\end{proof}

\begin{corollary}\label{cor:main}
$E$ belongs to the closed linear span of $\{G_k : k \geq 2\}$ in 
$H^2(\mathbb{C}_\alpha)$, for every $\alpha > 1$.
\end{corollary}

The membership of $E$ in the closed linear span of $\{G_k : k \geq 2\}$ also has
direct consequences for the location of the zeros of $\zeta$.

\begin{proposition}\label{main-prop}
For any $\alpha \geq 1/2$,
\[
E \in \overline{\mathrm{span}}\{G_k : k \geq 2\} \ \text{in } H^2(\mathbb{C}_\alpha) 
\implies \zeta(s) \neq 0, \quad \forall s \in \mathbb{C}_\alpha.
\]
\end{proposition}

\begin{proof}
Let $\rho$ be a zero of $\zeta$ in $\mathbb{C}_\alpha$. Then $G_k(\rho) = 0$ 
for every $k \geq 2$. Since $H^2(\mathbb{C}_\alpha)$ is a reproducing kernel
Hilbert space, convergence in norm implies pointwise convergence, so
$E(\rho) = 1/\rho = 0$, which is impossible.
\end{proof}

\subsection{A lower bound from the zeros of $\zeta$}\label{subsec:lower}

Proposition~\ref{main-prop} is qualitative: it turns convergence into a zero-free
half-plane. The reproducing kernel of $H^2(\mathbb{C}_\alpha)$ yields a quantitative
companion in the opposite direction, in which a single zero of $\zeta$ inside
$\mathbb{C}_\alpha$ produces an explicit, $n$-independent floor for the
approximation error. For $n \geq 2$ and $\alpha > 0$, define the
\emph{optimal approximation error}
\[
d_n(\alpha)^2 := \inf_{c_2, \dots, c_n \in \mathbb{C}}
\Bigl\| E - \sum_{k=2}^{n} c_k G_k \Bigr\|^2_\alpha,
\]
the best approximation of $E$ by the first $n-1$ generators. The restriction
$\alpha > 0$ is precisely what is needed for $E$ and the $G_k$ to belong to
$H^2(\mathbb{C}_\alpha)$: they are Mellin transforms of bounded functions supported
in $(0,1]$, and therefore belong to $H^2(\mathbb{C}_\alpha)$ for every $\alpha > 0$
(see the proof of Proposition~\ref{prop:exact}), whereas for $\alpha \leq 0$ they have
a pole at $s = 0$. Taking
$c_k = \mu(k)$ shows $d_n(\alpha)^2 \leq \|S_n - E\|^2_\alpha$, so any lower bound
for $d_n(\alpha)^2$ is a fortiori one for the Möbius error.

\begin{theorem}\label{thm:lower}
Let $\alpha > 0$ and suppose $\zeta$ has a zero $\rho$ with
$\mathrm{Re}(\rho) > \alpha$. Then for every $n \geq 2$ and all coefficients
$c_2, \dots, c_n \in \mathbb{C}$,
\[
\Bigl\| E - \sum_{k=2}^{n} c_k G_k \Bigr\|^2_\alpha
\ \geq\ \frac{|E(\rho)|^2}{\|k_\rho\|^2_\alpha}
\ =\ \frac{2\,(\mathrm{Re}\,\rho - \alpha)}{|\rho|^2}\ >\ 0.
\]
More generally, if $\rho_1, \dots, \rho_m$ are distinct zeros of $\zeta$ in
$\mathbb{C}_\alpha$, then for every $n \geq 2$,
\[
d_n(\alpha)^2 \ \geq\ v^{*} K^{-1} v, \qquad
K = \Bigl( \tfrac{1}{\rho_i + \overline{\rho_j} - 2\alpha} \Bigr)_{i,j=1}^{m}, \qquad
v = \bigl( 1/\rho_i \bigr)_{i=1}^{m}.
\]
In particular $d_n(\alpha)^2$ is bounded below by a positive constant independent
of $n$, so $d_n(\alpha)^2 \not\to 0$.
\end{theorem}

\begin{proof}
Since $\mathrm{Re}(\rho) > \alpha$ we have $\rho \in \mathbb{C}_\alpha$, and
$\zeta(\rho) = 0$ gives $G_k(\rho) = (k^{-\rho} - k^{-1})\zeta(\rho)/\rho = 0$ for
every $k$. Hence $F := E - \sum_{k=2}^{n} c_k G_k$ satisfies
$F(\rho) = E(\rho) = 1/\rho$. By the reproducing property
$F(\rho) = \langle F, k_\rho \rangle_\alpha$, so the Cauchy--Schwarz inequality gives
\[
\|F\|^2_\alpha \ \geq\ \frac{|\langle F, k_\rho \rangle_\alpha|^2}{\|k_\rho\|^2_\alpha}
= \frac{|E(\rho)|^2}{\|k_\rho\|^2_\alpha}.
\]
With $\|k_\rho\|^2_\alpha = k_\rho(\rho) = (2\,\mathrm{Re}\,\rho - 2\alpha)^{-1}$ and
$|E(\rho)|^2 = |\rho|^{-2}$, this is the first inequality.

For the general case, let $N = \{ g \in H^2(\mathbb{C}_\alpha) : g(\rho_i) = 0,\
1 \leq i \leq m \}$, whose orthogonal complement is
$V = \mathrm{span}\{ k_{\rho_1}, \dots, k_{\rho_m} \}$. Every combination
$\sum_k c_k G_k$ lies in $N$, since $G_k(\rho_i) = 0$; thus, with $P$ the orthogonal
projection onto $V$,
\[
\Bigl\| E - \sum_k c_k G_k \Bigr\|_\alpha
\ \geq\ \Bigl\| P\bigl( E - \textstyle\sum_k c_k G_k \bigr) \Bigr\|_\alpha
= \| P E \|_\alpha.
\]
A standard computation with the kernels gives
$\| P E \|^2_\alpha = v^{*} K^{-1} v$, where
$K_{ij} = \langle k_{\rho_j}, k_{\rho_i} \rangle_\alpha = k_{\rho_j}(\rho_i)
= (\rho_i + \overline{\rho_j} - 2\alpha)^{-1}$ and
$v_i = \langle E, k_{\rho_i} \rangle_\alpha = E(\rho_i) = 1/\rho_i$. Taking the
infimum over $c$ yields the bound on $d_n(\alpha)^2$; the right-hand side is a fixed
positive constant, so $d_n(\alpha)^2 \not\to 0$.
\end{proof}

\begin{corollary}\label{cor:sensor}
Let $\Theta = \sup\{\mathrm{Re}\,\rho : \rho \text{ a non-trivial zero of } \zeta\}$ and $\alpha_c = \inf\bigl\{\alpha \geq \tfrac12 : d_n(\alpha)^2 \to 0\bigr\}$. Then $\tfrac12 \leq \Theta \leq \alpha_c \leq 1$, and $\alpha_c = \tfrac12$ if and only if RH holds.
\end{corollary}

\begin{proof}
By Theorem~\ref{main-teo} the set $\{\alpha \geq \tfrac12 : d_n(\alpha)^2 \to 0\}$ is
nonempty and bounded below by $\tfrac12$, so $\alpha_c$ is well-defined and $\alpha_c \leq 1$.
For $\alpha < \Theta$ pick a zero $\rho$ with $\mathrm{Re}\,\rho > \alpha$;
Theorem~\ref{thm:lower} gives $d_n(\alpha)^2 \not\to 0$, whence $\Theta \leq \alpha_c$.
As $\Theta \geq \tfrac12$ is classical, $\tfrac12 \leq \Theta \leq \alpha_c \leq 1$. If
$\alpha_c = \tfrac12$ then $\Theta = \tfrac12$, which is RH. Conversely, under RH, Theorem~\ref{thm:bagchi} places $E$ in the closed linear span of
$\{G_k\}_{k\ge2}$ in $H^2(\mathbb{C}_{1/2})$. Since $d_n(1/2)$ is non-increasing
in $n$, this forces $d_n(1/2)^2\to0$; hence $\alpha=\tfrac12$ belongs to the set
$\{\alpha\ge\tfrac12:d_n(\alpha)^2\to0\}$, so that $\alpha_c\le\tfrac12$.
Combined with $\Theta\ge\tfrac12$, this gives $\alpha_c=\tfrac12$.
\end{proof}

\section{Numerical experiments}\label{sec:numerics}

The range $1/2 \leq \alpha \leq 1$ is not covered by Theorem~\ref{main-teo}; under
RH, convergence there follows from the estimates of Balazard and
Saias~\cite{Balazard} (see the end of Subsection~\ref{subsec:strip}), but only with a
weak bound. This section probes this range numerically. Two
questions guide the experiments: how the M\"obius error
$s_n(\alpha) = \|S_n - E\|^2_\alpha$ behaves as $n \to \infty$ across the strip, and
whether the obstruction on the critical line traced by Corollary~\ref{cor:sensor} is
intrinsic to the subspaces $\mathrm{span}\{G_2, \dots, G_n\}$ or only to the fixed
choice of coefficients $\mu(k)$. The first question is addressed in
Subsection~\ref{subsec:strip} and the second in Subsection~\ref{subsec:dn}. Both rest
on the computation of norms in $H^2(\mathbb{C}_\alpha)$, which we describe first.

\subsection{Computation of the norms}\label{subsec:method}

A direct evaluation of the boundary integral
\[
s_n(\alpha) = \frac{1}{2\pi}\int_{-\infty}^{\infty} \bigl|S_n(\alpha+it) - E(\alpha+it)\bigr|^2\,dt
\]
by numerical quadrature is costly and unreliable for large $n$: every evaluation of
the integrand involves $\zeta(\alpha+it)$ and the $n-1$ oscillating factors
$k^{-\alpha-it}$, and the integrand oscillates on scales that shrink as $n$ and $|t|$
grow. We avoid quadrature altogether by means of the following exact representation
of $s_n(\alpha)$ as a series with elementary terms.

\begin{proposition}\label{prop:exact}
Let $\alpha > 0$ and $n \geq 2$. For $m \geq 1$ put
\[
R_n(m) = 1 + \sum_{k=2}^{n} \mu(k)\,\frac{m \bmod k}{k}
\qquad\text{and}\qquad
w_m(\alpha) = \int_m^{m+1} \frac{du}{u^{2\alpha+1}} = \frac{m^{-2\alpha} - (m+1)^{-2\alpha}}{2\alpha}.
\]
Then
\begin{equation}\label{eq:exact}
s_n(\alpha) = \sum_{m=1}^{\infty} w_m(\alpha)\, R_n(m)^2 .
\end{equation}
\end{proposition}

\begin{proof}
Write $\{y\}$ for the fractional part of $y$ and, for $k \geq 2$, let
$\rho_k(x) = \{1/(kx)\} - \{1/x\}/k$ for $x > 0$. If $x > 1$, both $1/(kx)$ and $1/x$
lie in $(0,1)$, so $\rho_k(x) = 0$; moreover $|\rho_k| \leq 2$. For
$0 < \mathrm{Re}(s) < 1$ we have the classical formula~\cite[\S2.1]{titchmarsh},
written here after the change of variables $x \mapsto 1/x$,
\[
\int_0^\infty \{1/x\}\, x^{s-1}\, dx = -\frac{\zeta(s)}{s},
\]
and the substitution $y = kx$ gives
$\int_0^\infty \{1/(kx)\}\, x^{s-1}\,dx = -k^{-s}\zeta(s)/s$. Subtracting the two
identities,
\[
\int_0^1 \rho_k(x)\, x^{s-1}\, dx = -\bigl(k^{-s} - k^{-1}\bigr)\frac{\zeta(s)}{s} = -G_k(s),
\qquad 0 < \mathrm{Re}(s) < 1 .
\]
Since $\rho_k$ is bounded on $(0,1)$, the left-hand side converges absolutely and is
analytic on $\mathrm{Re}(s) > 0$, so by analytic continuation the identity holds on the
whole half-plane $\mathbb{C}_0$. Together with $E(s) = \int_0^1 x^{s-1}\,dx$ this gives
\[
E(s) - S_n(s) = \int_0^1 f_n(x)\, x^{s-1}\, dx,
\qquad f_n = 1 + \sum_{k=2}^{n} \mu(k)\,\rho_k,
\qquad \mathrm{Re}(s) > 0 .
\]
Let $\sigma \geq \alpha$. On the line $\mathrm{Re}(s) = \sigma$ the substitution
$x = e^{-y}$ turns this integral into the Fourier transform of
$g_\sigma(y) = f_n(e^{-y})\,e^{-\sigma y}$, $y > 0$, which lies in
$L^1(0,\infty) \cap L^2(0,\infty)$ because $f_n$ is bounded and $\sigma > 0$. By
Plancherel's theorem~\cite[Theorem~9.13]{rudin},
\[
\frac{1}{2\pi}\int_{-\infty}^{\infty} \bigl|(E - S_n)(\sigma+it)\bigr|^2\,dt
= \int_0^\infty |g_\sigma(y)|^2\,dy
= \int_0^1 |f_n(x)|^2\, x^{2\sigma-1}\,dx .
\]
Since $0 < x < 1$, the right-hand side is a non-increasing function of $\sigma$, and it
tends to its value at $\sigma = \alpha$ as $\sigma \to \alpha^+$ by monotone convergence.
Hence $E - S_n \in H^2(\mathbb{C}_\alpha)$ and, by the definition of the norm,
\[
s_n(\alpha) = \sup_{\sigma > \alpha} \int_0^1 |f_n(x)|^2\, x^{2\sigma-1}\,dx
= \int_0^1 |f_n(x)|^2\, x^{2\alpha-1}\,dx .
\]
Finally put $x = 1/u$, so that the last integral equals
$\int_1^\infty |f_n(1/u)|^2\, u^{-2\alpha-1}\,du$. For $u \in [m, m+1)$ the interval
$(m, u]$ contains no integer, hence no multiple of $k$; thus
$\{u/k\} = u/k - \lfloor m/k \rfloor$ and $\{u\} = u - m$, which give
$\rho_k(1/u) = m/k - \lfloor m/k \rfloor = (m \bmod k)/k$ and $f_n(1/u) = R_n(m)$.
Splitting the integral over the intervals $[m, m+1)$ yields~\eqref{eq:exact}.
\end{proof}

\begin{remark}\label{rem:exact}
(i) For $\alpha = 1/2$ the weights are $w_m(1/2) = 1/(m(m+1))$, and
Proposition~\ref{prop:exact} turns the approximation of $E$ in
$H^2(\mathbb{C}_{1/2})$ into an approximation problem for the constant sequence
$\mathbf{1}$ by the sequences $r_k/k$, $r_k(m) = m \bmod k$, in a weighted
$\ell^2$-space; compare with the setting of Balazard and Saias~\cite{Balazard}
recalled in the introduction. The proof also shows that
$s_n(1/2) = \|f_n\|^2_{L^2(0,1)}$, which makes explicit the link with the
setting of B\'aez-Duarte. Indeed, put $\chi = \chi_{(0,1]}$,
$M_1(n) = \sum_{k \leq n} \mu(k)/k$, and extend $f_n$ by $0$ to $(1,\infty)$.
Since $\sum_{k=2}^{n} \mu(k)/k = M_1(n) - 1$, the definition of the $\rho_k$ gives
\[
f_n(x) = \chi(x) + \sum_{j=1}^{n} \mu(j) \Bigl\{ \frac{1}{jx} \Bigr\}
- M_1(n) \Bigl\{ \frac{1}{x} \Bigr\}, \qquad x > 0,
\]
both sides vanishing for $x > 1$, where the right-hand side equals
$M_1(n)/x - M_1(n)/x$. As $\{1/x\} \in L^2(0,\infty)$ and $M_1(n) \to 0$, the
convergence $s_n(1/2) \to 0$ would force the partial sums of~\eqref{eq:bd-series} to
converge to $\chi$ in $L^2(0,\infty)$, contradicting
\cite[Theorem~2.2]{baez-unitary}. Hence $s_n(1/2) \not\to 0$; and since, by
Theorem~\ref{main-teo}, any limit of $(S_n)$ in $H^2(\mathbb{C}_{1/2})$ would
coincide with $E$ on $\mathbb{C}_1$, hence on $\mathbb{C}_{1/2}$, the sequence
$(S_n)$ does not converge in $H^2(\mathbb{C}_{1/2})$.

(ii) Since $\sum_{k \leq m} \mu(k) \lfloor m/k \rfloor = 1$ for every
$m \geq 1$~\cite[Theorem~3.12]{apostol}, and $\lfloor m/k \rfloor = 0$ for $k > m$, for
$m \leq n$ we get
\[
R_n(m) = 1 + \sum_{k=1}^{n} \mu(k)\Bigl(\frac{m}{k} - \Bigl\lfloor \frac{m}{k} \Bigr\rfloor\Bigr)
= m\, M_1(n), \qquad M_1(n) = \sum_{k \leq n} \frac{\mu(k)}{k},
\]
and therefore
\begin{equation}\label{eq:lowerM1}
s_n(\alpha) \geq M_1(n)^2 \sum_{m=1}^{n} m^2\, w_m(\alpha),
\qquad
\sum_{m=1}^{n} m^2\, w_m(\alpha) \sim \frac{n^{2-2\alpha}}{2-2\alpha}
\quad \Bigl(\tfrac12 \leq \alpha < 1\Bigr).
\end{equation}
Thus the decay of $s_n(\alpha)$ is limited from below by the size of the partial sums
$M_1(n)$ of $\sum_k \mu(k)/k$. Under RH one has $M_0(x)=\sum_{k\le x}\mu(k)=O(x^{1/2+\varepsilon})$
\cite[Theorem~14.25(C)]{titchmarsh}; since $\sum_{k}\mu(k)/k=0$, partial
summation gives $M_1(x)=M_0(x)/x-\int_x^\infty M_0(t)\,t^{-2}\,dt$, whence
$M_1(n)=O(n^{-1/2+\varepsilon})$; if $M_1(n)$ is of exact order $n^{-1/2}$ along a
sequence of values of $n$, the right-hand side of~\eqref{eq:lowerM1} is of order
$n^{1-2\alpha}$ along that sequence.
\end{remark}

\paragraph{Implementation.}
We evaluate the series~\eqref{eq:exact} in three parts. The terms with $m \leq m_0$ are
summed directly. On $m_0 < m \leq M$ the terms are grouped into blocks $[a, b)$ of
relative length $(b-a)/a \leq h$, and the weight is expanded to second order about the
centre $c$ of each block,
\[
\sum_{a \leq m < b} w_m(\alpha)\, R_n(m)^2 \approx
\sum_{j=0}^{2} \frac{w^{(j)}(c)}{j!} \sum_{a \leq m < b} (m - c)^j R_n(m)^2,
\qquad w(x) = \frac{x^{-2\alpha} - (x+1)^{-2\alpha}}{2\alpha},
\]
which has relative error $O(h^3)$ and reduces the dependence on $\alpha$ to three
moments of $R_n^2$ per block, shared by all values of $\alpha$. For $m > M$ we replace
$R_n(m)^2$ by its mean value $\overline{R_n^2}$ over a period of $m \mapsto R_n(m)$, so
that the tail contributes $\overline{R_n^2}\,(M+1)^{-2\alpha}/(2\alpha)$. This
replacement is not exact, since the period $\mathrm{lcm}(2, \dots, n)$ of
$m \mapsto R_n(m)$ is far larger than $M$; its effect is part of the truncation
error estimated below. The mean itself is known exactly. Let $g = \gcd(j,k)$. For $m$ uniformly distributed
modulo $\mathrm{lcm}(j,k)$, the Chinese remainder theorem gives
$m \bmod j = ga + r$ and $m \bmod k = gb + r$, where $r = m \bmod g$,
$a \in \{0, \dots, j/g - 1\}$ and $b \in \{0, \dots, k/g - 1\}$ are independent and
uniformly distributed; hence the covariance of $m \bmod j$ and $m \bmod k$ equals
the variance of $r$, namely $(g^2 - 1)/12$. Writing
$\gcd(j,k)^2 = \sum_{d \mid \gcd(j,k)} J_2(d)$, where
$J_2(d) = d^2 \prod_{p \mid d} (1 - p^{-2})$ is Jordan's totient, characterised by
$\sum_{d \mid q} J_2(d) = q^2$ for $q \geq 1$ \cite[Ch.~2, Exercise~17]{apostol}
(the term $d = 1$ cancels the $-1$), one finds
\[
\overline{R_n^2} = \Bigl(1 + \sum_{k=2}^{n} \mu(k)\,\frac{k-1}{2k}\Bigr)^2
+ \frac{1}{12} \sum_{d \geq 2} J_2(d) \Bigl(\sum_{\substack{2 \leq k \leq n \\ d \mid k}} \frac{\mu(k)}{k}\Bigr)^2 .
\]
Since $R_n(m) = R_{n-1}(m) + \mu(n)\,(m \bmod n)/n$, passing from $n-1$ to $n$ costs a
single pass over $m \leq M$, and only the $n$ with $\mu(n) \neq 0$ contribute. All the
values $s_n(\alpha)$ with $n \leq N$ are therefore obtained with $O(NM)$ elementary
operations, without evaluating $\zeta$ and without numerical quadrature. Since all
terms of~\eqref{eq:exact} are non-negative, no cancellation occurs when the series is
summed. All computations are carried out in IEEE double precision, without
compensated summation; two implementations of the recursion that accumulate $R_n(m)$
in different orders agree to within $5 \cdot 10^{-9}$ (relative), far below the
truncation error estimated below.
The computation was organised so that each block of values of $m$ stays in cache while it
runs through all $n$, the blocks being processed in parallel. A Fortran~2008 program
with OpenMP computed $s_n(\alpha)$ for the $51$ values
$\alpha = 0.50, 0.51, \dots, 1.00$ and the $3001$ logarithmically spaced integers
$n_j$, $0 \leq j \leq 3000$, where $n_j$ is the integer nearest to $10^{3 + j/1000}$
(these integers are distinct, and $n_{3000} = 10^6$), in $40$ minutes on the machine
described at the end of this subsection,
with $m_0 = 10^4$, $h = 10^{-3}$ and $M = 3 \cdot 10^8$.
The block expansion changes the results by less than $10^{-9}$ (relative) with respect
to direct summation. To control the truncation error we repeated the whole
computation with $M = 10^8$ and with $M = 10^9$. With respect to the run with
$M = 10^9$, the values obtained with $M = 3 \cdot 10^8$ differ by at most $10^{-8}$
for $n < 10^4$, $7.4 \cdot 10^{-7}$ for $10^4 \leq n < 10^5$ and $4.8 \cdot 10^{-5}$
for $10^5 \leq n \leq 10^6$, while those obtained with $M = 10^8$ differ from the run
with $M = 3 \cdot 10^8$ by up to $3.5 \cdot 10^{-4}$. These differences are consistent
with a truncation error decreasing like $M^{-1.8}$, so that the relative error of the
values reported below, computed with $M = 3 \cdot 10^8$, is about $5 \cdot 10^{-5}$ at
most; the entries of Table~\ref{tab:envelope} change by at most $1.1 \cdot 10^{-5}$ when
$M = 10^9$ is used instead.
As an independent check, the values of $s_n(1/2)$ given by~\eqref{eq:exact} agree
with those obtained from the Gram data described below, which are computed from the
boundary integral and the zeta function, to within $3.4 \cdot 10^{-4}$ (relative) for
$2 \leq n \leq 5120$, the median difference being $6 \cdot 10^{-5}$. The same check on
the line $\mathrm{Re}(s) = 3/4$, with the Gram data computed there in the same way,
gives agreement with~\eqref{eq:exact} to within $10^{-4}$ for $2 \leq n \leq 5120$, the
median difference being $10^{-5}$.

\paragraph{Optimal error.} The optimal error of Subsection~\ref{subsec:lower} is
computed on the critical line as
$d_n(1/2)^2 = \|E\|^2_{1/2} - b^{\top} A^{-1} b$, with Gram data
$A = (\langle G_j, G_k \rangle_{1/2})_{j,k=2}^{n}$ and
$b = (\langle E, G_k \rangle_{1/2})_{k=2}^{n}$. Because $E$ and the $G_k$ are
real on the real axis, they satisfy the reflection symmetry
$G_k(\bar s) = \overline{G_k(s)}$ (and likewise for $E$), so each such inner
product reduces to a real half-line integral,
\[
\langle G_j, G_k \rangle_{1/2}
= \frac{1}{\pi} \int_0^\infty \mathrm{Re}\,\bigl[G_j(\tfrac12+it)\,\overline{G_k(\tfrac12+it)}\bigr]\,dt,
\]
and likewise for $\langle E, G_k \rangle_{1/2}$; hence $A$ and $b$ are real and the
minimising coefficients may be taken real. These boundary integrals are evaluated with
the trapezoidal rule on a fixed grid of $33\,600$ nodes, $8000$ of them equally spaced
on $[0, 80)$ and $25\,600$ geometrically spaced on $[80, 3 \cdot 10^6]$, with
$\zeta(\tfrac12 + it)$ computed by the Euler--Maclaurin formula with $N \geq 2|s|/\pi$
terms and twelve Bernoulli corrections, whose remainder is negligible in double
precision. On this grid $\|E\|^2_{1/2} = 1$ is reproduced to $1.1 \cdot 10^{-7}$, the
values of $\langle E, G_k \rangle_{1/2}$ agree to within $10^{-4}$ (relative) with the
exact value $G_k(1) = -(\log k)/k$, which holds because $E$ is the reproducing kernel of
$H^2(\mathbb{C}_{1/2})$ at $s = 1$, and the values of $s_n(1/2)$ computed from the Gram
data agree with~\eqref{eq:exact} as stated above; these checks control the accuracy
of the inner products. The Gram matrix for $n = 5120$ is
factorised as $A = U^{\top} U$ (Cholesky), and with $y = U^{-\top} b$ the leading blocks
of the factorisation give
\[
d_n(1/2)^2 = \|E\|^2_{1/2} - \sum_{i=1}^{n-1} y_i^2, \qquad 2 \leq n \leq 5120,
\]
for all $n$ at once. Since $A$ is ill-conditioned, the accuracy of the inner
products does not by itself bound the error of $d_n(1/2)^2$. We therefore also
computed $d_n(1/2)^2$ for $n \leq 1280$ from the exact Gram data that follow from the
proof of Proposition~\ref{prop:exact},
\[
\langle G_j, G_k \rangle_{1/2} = \frac{1}{jk} \sum_{m=1}^{\infty}
\frac{(m \bmod j)(m \bmod k)}{m(m+1)},
\qquad
\langle E, G_k \rangle_{1/2} = -\frac{\log k}{k},
\]
where the series is summed exactly over one period $L = \mathrm{lcm}(j,k)$ by means of
$\sum_{q \geq 0} \bigl((qL + r)(qL + r + 1)\bigr)^{-1} = \bigl(\psi((r+1)/L) - \psi(r/L)\bigr)/L$,
$\psi$ being the digamma function; closed formulas for these inner products were
obtained by Vasyunin~\cite{vasyunin}. The two computations of $d_n(1/2)^2$ agree to
within $2 \cdot 10^{-4}$ (relative) for $n < 640$ and to within $5 \cdot 10^{-4}$ for
$640 \leq n \leq 1280$, the difference growing slowly with $n$ as $A$ becomes more
ill-conditioned; for larger $n$ the comparison with~\cite{BBLS} in
Subsection~\ref{subsec:dn} provides a further check. The exact Gram data also give
$s_n(1/2)$, which agrees with the implementation of~\eqref{eq:exact} described above to
within $4 \cdot 10^{-8}$.

\paragraph{Computing environment.}
All the computations of this section were carried out with programs written in
Fortran~2008 and parallelised with OpenMP, on a desktop computer with an Intel Core
i7-13700KF processor (16 cores, 8 performance and 8 efficiency cores, 24 hardware
threads, maximum clock frequency 5.4~GHz, 30~MB of L3 cache) and 32~GB of RAM, running
Ubuntu~26.04~LTS (Linux kernel~7.0). The programs were compiled with GNU Fortran~15.2
with the options \texttt{-O3 -march=native -fopenmp} and run with $20$ OpenMP threads;
the matrix products and the Cholesky factorisation of the optimal-error computation use
the BLAS and LAPACK libraries. The computation of $d_n(1/2)^2$ for $2 \leq n \leq 5120$
took less than two minutes. The figures were produced with Python (Matplotlib) from the
output of the Fortran programs. The code generating all figures and tables of this
section is available from the authors upon request.

\subsection{The M\"obius error across the strip}\label{subsec:strip}

\begin{figure}[htb!]
\centering
\includegraphics[width=0.9\textwidth]{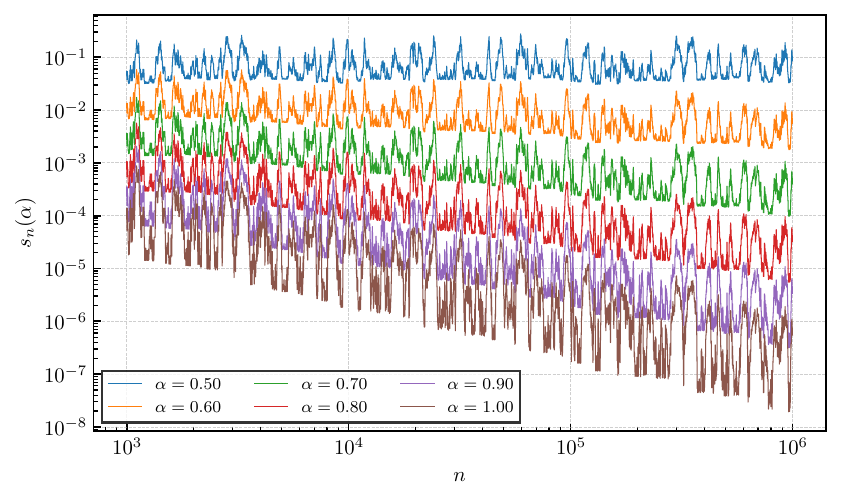}
\caption{The M\"obius error $s_n(\alpha)$ for $\alpha = 0.50, 0.60, \dots, 1.00$ and
$10^3 \leq n \leq 10^6$, in log--log scale.}
\label{fig:sn}
\end{figure}

\begin{figure}[htb!]
\centering
\includegraphics[width=0.7\textwidth]{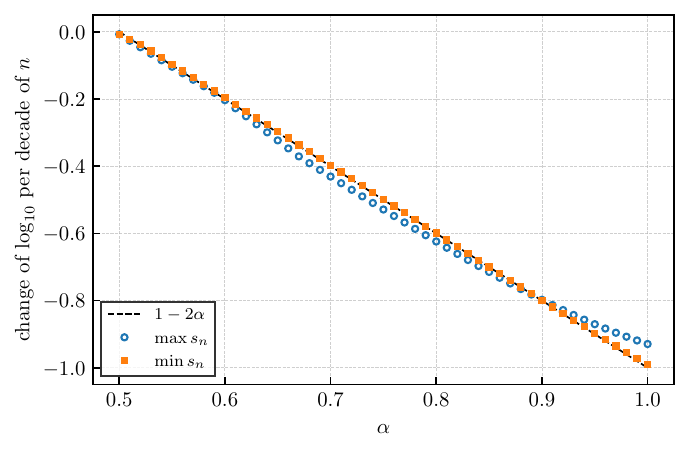}
\caption{Change of $\log_{10} \max s_n(\alpha)$ and of $\log_{10} \min s_n(\alpha)$ per
decade of $n$, from $10^3 \leq n < 10^4$ to $10^5 \leq n \leq 10^6$, for
$\alpha = 0.50, 0.51, \dots, 1.00$, compared with $1 - 2\alpha$ (dashed).}
\label{fig:exponent}
\end{figure}

We computed $s_n(\alpha)$ for $\alpha = 0.50, 0.51, \dots, 1.00$ and
$10^3 \leq n \leq 10^6$, at $1000$ logarithmically spaced values of $n$ per decade.
Figure~\ref{fig:sn} shows the sequences for $\alpha = 0.50, 0.60, \dots, 1.00$. They
oscillate strongly with $n$: within each decade the largest value exceeds the smallest
by a factor of about $8$ on the critical line, and by a factor of several hundred for
$\alpha = 1$. The oscillations occur at the same values of $n$ for all $\alpha$, and
the value of $s_n(\alpha)$ at an isolated $n$ is therefore not representative of its
behaviour. We summarise each
sequence instead by its largest and smallest values over each decade of $n$
(Table~\ref{tab:envelope}), and measure its decay by the change of
$\log_{10} \max s_n$ and $\log_{10} \min s_n$ per decade between the first and the last
decade (Figure~\ref{fig:exponent}). Since it is based on extreme values, this
two-point estimate is only indicative.

\begin{table}[htb!]
\centering
\footnotesize
\begin{tabular}{|c|cc|cc|cc|}
\hline
\multirow{2}{*}{$\alpha$} & \multicolumn{2}{c|}{$10^{3} \leq n < 10^{4}$} & \multicolumn{2}{c|}{$10^{4} \leq n < 10^{5}$} & \multicolumn{2}{c|}{$10^{5} \leq n \leq 10^{6}$} \\
\cline{2-7}
 & $\max s_n$ & $\min s_n$ & $\max s_n$ & $\min s_n$ & $\max s_n$ & $\min s_n$ \\
\hline\hline
0.50 & $2.617\mathrm{e}{-1}$ & $3.181\mathrm{e}{-2}$ & $2.772\mathrm{e}{-1}$ & $3.408\mathrm{e}{-2}$ & $2.530\mathrm{e}{-1}$ & $3.079\mathrm{e}{-2}$ \\
0.55 & $1.219\mathrm{e}{-1}$ & $1.189\mathrm{e}{-2}$ & $9.722\mathrm{e}{-2}$ & $1.086\mathrm{e}{-2}$ & $7.554\mathrm{e}{-2}$ & $7.605\mathrm{e}{-3}$ \\
0.60 & $5.866\mathrm{e}{-2}$ & $4.390\mathrm{e}{-3}$ & $3.954\mathrm{e}{-2}$ & $3.213\mathrm{e}{-3}$ & $2.298\mathrm{e}{-2}$ & $1.776\mathrm{e}{-3}$ \\
0.65 & $3.162\mathrm{e}{-2}$ & $1.632\mathrm{e}{-3}$ & $1.683\mathrm{e}{-2}$ & $9.559\mathrm{e}{-4}$ & $7.140\mathrm{e}{-3}$ & $4.171\mathrm{e}{-4}$ \\
0.70 & $1.740\mathrm{e}{-2}$ & $6.106\mathrm{e}{-4}$ & $7.338\mathrm{e}{-3}$ & $2.857\mathrm{e}{-4}$ & $2.393\mathrm{e}{-3}$ & $9.806\mathrm{e}{-5}$ \\
0.75 & $9.798\mathrm{e}{-3}$ & $2.296\mathrm{e}{-4}$ & $3.291\mathrm{e}{-3}$ & $8.578\mathrm{e}{-5}$ & $8.578\mathrm{e}{-4}$ & $2.316\mathrm{e}{-5}$ \\
0.80 & $5.662\mathrm{e}{-3}$ & $8.677\mathrm{e}{-5}$ & $1.528\mathrm{e}{-3}$ & $2.584\mathrm{e}{-5}$ & $3.197\mathrm{e}{-4}$ & $5.499\mathrm{e}{-6}$ \\
0.85 & $3.370\mathrm{e}{-3}$ & $3.293\mathrm{e}{-5}$ & $7.396\mathrm{e}{-4}$ & $7.811\mathrm{e}{-6}$ & $1.254\mathrm{e}{-4}$ & $1.313\mathrm{e}{-6}$ \\
0.90 & $2.073\mathrm{e}{-3}$ & $1.255\mathrm{e}{-5}$ & $3.768\mathrm{e}{-4}$ & $2.368\mathrm{e}{-6}$ & $5.263\mathrm{e}{-5}$ & $3.157\mathrm{e}{-7}$ \\
0.95 & $1.323\mathrm{e}{-3}$ & $4.793\mathrm{e}{-6}$ & $2.039\mathrm{e}{-4}$ & $7.196\mathrm{e}{-7}$ & $2.405\mathrm{e}{-5}$ & $7.693\mathrm{e}{-8}$ \\
1.00 & $8.784\mathrm{e}{-4}$ & $1.838\mathrm{e}{-6}$ & $1.181\mathrm{e}{-4}$ & $2.192\mathrm{e}{-7}$ & $1.217\mathrm{e}{-5}$ & $1.926\mathrm{e}{-8}$ \\
\hline
\end{tabular}
\caption{Largest and smallest values of the M\"obius error $s_n(\alpha)$ over
each decade of $n$, for $\alpha = 0.50, 0.55, \dots, 1.00$, computed at $3001$
logarithmically spaced values of $n \in [10^3, 10^6]$. The maxima and minima are
taken over these sampled values of $n$.}
\label{tab:envelope}
\end{table}

On the critical line the M\"obius error does not decay: over the three decades the
maxima of $s_n(1/2)$ are $0.262$, $0.277$ and $0.253$, and the minima
$3.18 \cdot 10^{-2}$, $3.41 \cdot 10^{-2}$ and $3.08 \cdot 10^{-2}$. This is
consistent with the result of B\'aez-Duarte~\cite{baez-unitary}, by which
$s_n(1/2) \not\to 0$ (Remark~\ref{rem:exact}). For every $\alpha > 1/2$, on the other hand, both the maxima and
the minima decrease geometrically from one decade to the next, and
Figure~\ref{fig:exponent} shows that the rate is very close to $10^{1-2\alpha}$ per
decade: over the whole range the observed changes of $\log_{10}$ per decade differ from
$1 - 2\alpha$ by less than $0.08$. The data are thus consistent with
\[
s_n(\alpha) \asymp n^{1-2\alpha}, \qquad \tfrac12 \leq \alpha < 1,
\]
for the envelope of the oscillations (for $\alpha = 1$ see below), and in particular with $s_n(\alpha) \to 0$ for
every $\alpha > 1/2$. The convergence is very slow near the critical line: for
$\alpha = 0.55$ the error decreases only by a factor $10^{0.1} \approx 1.26$ per decade,
so that for moderate $n$ the behaviour for $\alpha$ slightly above $1/2$ is hard to tell
apart from that on the critical line. The exponent $1 - 2\alpha$ is precisely the one
given by the lower bound~\eqref{eq:lowerM1} when $M_1(n)$ is of order $n^{-1/2}$, and
the largest peaks of $s_n(\alpha)$ in Figure~\ref{fig:sn} occur where $n M_1(n)^2$ is
large. Near $\alpha = 1$ the maxima decay more slowly than $n^{1-2\alpha}$
(Figure~\ref{fig:exponent}). This is consistent with~\eqref{eq:lowerM1}: for
$\alpha = 1$ the sum $\sum_{m \leq n} m^2 w_m(1)$ grows like $\log n$, so that along
a sequence with $M_1(n)$ of exact order $n^{-1/2}$ the lower bound is of order
$n^{-1}\log n$, whose slope per decade over the computed range is about $-0.91$,
close to the observed $-0.93$.

These observations must be read with care. By Proposition~\ref{main-prop}, a proof that
$s_n(\alpha) \to 0$ for some $\alpha < 1$ would give a zero-free half-plane
$\mathbb{C}_\alpha$, far beyond the known zero-free regions of $\zeta$. The behaviour
observed is the one expected under RH. Indeed, by Remark~\ref{rem:exact},
$f_n = 1 + S_n^{*} - M_1(n)\{1/x\}$ on $(0,1]$ with
$S_n^{*}(x) = \sum_{k \leq n} \mu(k)\{1/(kx)\}$, and if RH holds the estimates of
Balazard and Saias~\cite[Lemmes~4 et~5]{Balazard} give
$\int_0^1 |1 + S_n^{*}(x)|^2 x^{2\alpha-1}\,dx \ll_\alpha n^{-(\alpha-1/2)/3}$ for
$1/2 < \alpha \leq 1$; together with $M_1(n) = O(n^{-1/2+\varepsilon})$, this yields
$s_n(\alpha) = O(n^{-(\alpha-1/2)/3})$, and in particular $s_n(\alpha) \to 0$, for
every $\alpha > 1/2$. The decay like $n^{1-2\alpha}$ observed above is much faster,
and matches the order of the
lower bound~\eqref{eq:lowerM1}. Nevertheless, no finite computation can rule out the
alternative: by Theorem~\ref{thm:lower}, a
zero $\rho$ with $\mathrm{Re}\,\rho > \alpha$ only forces
$s_n(\alpha) \geq 2(\mathrm{Re}\,\rho - \alpha)/|\rho|^2$, and since every zero with
$|\mathrm{Im}\,\rho| \leq 3 \cdot 10^{12}$ lies on the critical line~\cite{platt}, this
floor would be smaller than $10^{-24}$.

\subsection{The optimal error on the critical line}\label{subsec:dn}

The critical line is the one place in the range studied where $S_n$ demonstrably
fails: by B\'aez-Duarte~\cite{baez-unitary} the Möbius coefficients do not sum to
$E$ at $\alpha = 1/2$ (Remark~\ref{rem:exact}), and the undamped oscillations of Figure~\ref{fig:sn}
are the price of this fixed choice. Since it is exactly at $\alpha = 1/2$ that
Theorem~\ref{thm:bagchi} formulates the RH, the natural question is what this
failure rules out. The answer is: nothing about the subspaces
$\mathrm{span}\{G_2, \dots, G_n\}$ themselves, since $\mu(k)$ is only one
admissible choice of coefficients, and it is the subspaces, not the choice, that
Theorem~\ref{thm:lower} ties to the zeros of $\zeta$. This is what motivates the
present experiment: we replace $s_n(1/2)$ by the optimal error $d_n(1/2)^2$ of
Subsection~\ref{subsec:lower}, the true distance from $E$ to
$\mathrm{span}\{G_2, \dots, G_n\}$ on the critical line, which removes the
coefficients from the picture by optimising over all of them. It is computed by
the Gram matrix method of Subsection~\ref{subsec:method}.

\begin{figure}[htb!]
\centering
\begin{minipage}[c]{0.54\textwidth}
\centering
\includegraphics[width=\linewidth]{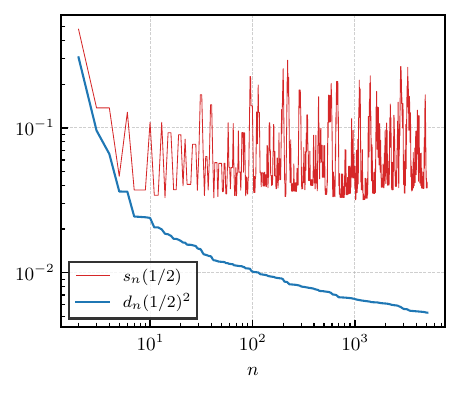}
\end{minipage}\hfill
\begin{minipage}[c]{0.42\textwidth}
\centering
\renewcommand{\arraystretch}{1.2}
\begin{tabular}{|r|c|c|}
\hline
$n$ & $d_n(1/2)^2$ & $s_n(1/2)$ \\
\hline
\hline
$40$   & $1.274\mathrm{e}{-2}$ & $1.442\mathrm{e}{-1}$ \\
$160$  & $9.357\mathrm{e}{-3}$ & $6.847\mathrm{e}{-2}$ \\
$640$  & $7.023\mathrm{e}{-3}$ & $3.659\mathrm{e}{-2}$ \\
$1280$ & $6.366\mathrm{e}{-3}$ & $3.474\mathrm{e}{-2}$ \\
$2560$ & $5.932\mathrm{e}{-3}$ & $4.314\mathrm{e}{-2}$ \\
$5120$ & $5.283\mathrm{e}{-3}$ & $3.851\mathrm{e}{-2}$ \\
\hline
\end{tabular}
\end{minipage}

\medskip
\begin{minipage}[t]{0.54\textwidth}
\caption{The optimal error $d_n(1/2)^2$ (blue) and the Möbius error
$s_n(1/2)$ (red) on the critical line, \mbox{$2 \leq n \leq 5120$}, in log--log scale.}
\label{fig:dn}
\end{minipage}\hfill
\begin{minipage}[t]{0.42\textwidth}
\captionof{table}{The same two errors at selected $n$; for $n \leq 1280$ the values
of $d_n(1/2)^2$ are those of the exact Gram data.}
\label{tab:dn}
\end{minipage}
\end{figure}

Figure~\ref{fig:dn} contrasts the two sequences, and Table~\ref{tab:dn} records
their values at selected $n$. While the M\"obius error $s_n(1/2)$ keeps oscillating
between about $3 \cdot 10^{-2}$ and $3 \cdot 10^{-1}$, the optimal error
$d_n(1/2)^2$, which is non-increasing in $n$ by definition, decreases slowly:
for $100 \leq n \leq 5120$, $d_n(1/2)^2 \log n$ stays between $0.044$ and $0.048$.
Under RH, $d_n(1/2)^2 \to 0$ by Theorem~\ref{thm:bagchi}, whereas
$s_n(1/2) \not\to 0$ unconditionally (Remark~\ref{rem:exact}); thus, under RH, the
failure of $S_n$ on the critical line is due to the fixed choice $c_k = \mu(k)$
and not to the subspaces $\mathrm{span}\{G_2, \dots, G_n\}$. Our computations
illustrate this contrast; they cannot, of course, decide whether
$d_n(1/2)^2 \to 0$ (see the end of Subsection~\ref{subsec:strip}).

The distance $d_n(1/2)$ has been studied before in the Nyman--Beurling setting. By
the proof of Proposition~\ref{prop:exact}, $d_n(1/2)$ is the distance in
$L^2(0,\infty)$ from $\chi = \chi_{(0,1]}$ to the functions
$\sum_{k=1}^{n} c_k \{1/(kx)\}$ with $\sum_{k} c_k/k = 0$. Dropping this constraint
gives the distance studied by B\'aez-Duarte, Balazard, Landreau and
Saias~\cite{BBLS}, which is in turn bounded below by their $D(1/n)$. Hence, by
\cite[Th\'eor\`eme]{BBLS} as sharpened by Burnol~\cite[Theorem~1.3]{burnol},
\[
\liminf_{n \to \infty} d_n(1/2)^2 \log n \ \geq\ 2 + \gamma - \log 4\pi = 0.0461\ldots
\]
unconditionally: if RH fails the bound is trivial, since then $d_n(1/2)$ does not tend
to $0$, while under RH Burnol's bound $\sum_\rho m_\rho^2/|\rho|^2$, where $m_\rho$ is
the multiplicity of the zero $\rho$, is at least
$\sum_\rho m_\rho/|\rho|^2 = 2 + \gamma - \log 4\pi$. Being a bound on the $\liminf$,
it is compatible with the values of $d_n(1/2)^2 \log n$ slightly below
$2 + \gamma - \log 4\pi$ found above for finite $n$. In~\cite[\S8]{BBLS} the unconstrained distance was computed for
$n \leq 10^4$ from Vasyunin's closed formulas for the inner
products~\cite{vasyunin}, which led them to conjecture
\cite[Conjecture~2]{BBLS} that $d_n^2 \sim (2 + \gamma - \log 4\pi)/\log n$. Our values,
obtained by a different method (the exact Gram data of Subsection~\ref{subsec:method}
for $n \leq 1280$ and quadrature on the critical line beyond), are consistent with
these results, although the unconstrained distance of~\cite{BBLS} does not exceed
ours: for $100 \leq n \leq 5120$, $d_n(1/2)^{-2}$ differs by at most
$5.1\%$ from the least-squares fit $21.99 \log n$ of~\cite[\S8]{BBLS}. This agreement also provides
an independent check of the computation described in
Subsection~\ref{subsec:method}.

\section*{Declaration of Generative AI Use}
During the preparation of this work the authors used Claude (Anthropic, model Claude
Opus~5.5) to improve the English language and readability of the manuscript, to assist
in writing and debugging the Fortran and Python code used for the numerical experiments
and figures, to check the proofs for errors, and to search the literature. It was also
used to shorten the abstract to the journal's word limit and to check the compilation of
the \LaTeX{} source. All mathematical results, proofs and numerical computations were
conceived, verified and validated by the authors, who reviewed and edited all AI-assisted
output and take full responsibility for the content of the manuscript.

\section*{Acknowledgements}
The authors thank the Vice-Rectorate for Research, Outreach and Social Projection and
the Applied Mathematics research group of the University of the Atlantic for providing
the facilities of the Center for Mathematical Modeling, where this work was carried
out.

\section*{Disclosure statement}
No potential conflict of interest was reported by the authors.

\section*{Notes on contributors}
\textit{J. Manzur} is a faculty member at the Department of Mathematics of the
Universidad del Atl\'antico, working in operator theory on Hilbert spaces of analytic
functions.

\medskip
\noindent
\textit{G. Quintero} is a full-time faculty member at the Department of Mathematics of
the Universidad del Atl\'antico. His expertise is in nonlinear optimization, with a
focus on the development and analysis of numerical optimization methods and
algorithms.

\FloatBarrier

\end{document}